\documentclass[leqno]{article}
\usepackage{amssymb, amsmath, amsfonts}

\newtheorem{theorem}{Theorem}[section]
\newtheorem{corollary}[theorem]{Corollary}
\newtheorem{lemma}[theorem]{Lemma}
\newtheorem{proposition}[theorem]{Proposition}
\newtheorem{remark}[theorem]{Remark}

\numberwithin{equation}{section}

\def\sqr#1#2{{\vcenter{\vbox{\hrule height.#2pt
    \hbox{\vrule width.#2pt height#1pt \kern#1pt
    \vrule width.#2pt}
    \hrule height.#2pt}}}}

\def\ga{\gamma}

\def\De{\Delta}
\def\de{\delta}

\def\la{\lambda}
\def\eps{\varepsilon}
\def\ka{\kappa}

\def\om{\omega}
\def\Om{\Omega}

\def\pa{\partial}

\def\bB{\bar{B}}

\def\bOm{\bar\Om}

\def\tV{\tilde V}

\def\tI{\tilde I}

\font\bbb=msbm10

\def\Z{\hbox{{\bbb Z}}}
\def\R{\hbox{{\bbb R}}}
\def\C{\hbox{{\bbb C}}}

\def\cH{{\cal H}}

\def\Ker{\, {\rm Ker}\, }
\def\capa{\, {\rm cap}\, }

\def\Lip{\, {\rm Lip}}

\def\Tr{\, {\rm Tr}\, }

\def\supp{\hbox{supp}\;}

\def\Int{{\rm Int}\,}
\def\cM{{\cal M}}

\def\dist{{\rm dist}}

\def\ms{\medskip}

\begin{document}
\title{Gauge optimization and spectral properties\\  
of magnetic Schr\"odinger operators}

\author{
{\bf Vladimir Kondratiev\thanks{Research partially supported by 
Matthews professorship fund of Northeastern University}}
\smallskip\\
Moscow State University, Moscow, Russia
\bigskip\\
{\bf Vladimir Maz'ya\thanks{Research partially supported by the Department of
Mathematics and the Robert G. Stone Fund at Northeastern
University}}\\
Ohio State University,
Columbus, OH, USA
\\
University of Liverpool, Liverpool,  UK
\\
Link\"oping University,
Sweden
\bigskip\\
{\bf Mikhail Shubin
\thanks{Research partially supported by NSF grant DMS-0107796} }
\\
Northeastern University,
Boston, MA,
USA
}

\date{}

\maketitle

\begin{abstract}

We establish new necessary and sufficient conditions
for the discreteness of spectrum 
and strict positivity
of magnetic Schr\"odinger operators with a positive scalar potential. 
They extend earlier results by  Maz'ya and Shubin (2005),
which were obtained in case when there is no magnetic field. 
We also derive two-sided estimates for the bottoms of spectrum and essential spectrum,
extending results by Maz'ya and Otelbaev (1977). 
The main idea is to optimize the gauge of the magnetic field, 
thus reducing the quadratic form 
to  one without magnetic field 
(but with an appropriately adjusted scalar potential). 
\end{abstract}

\section{Introduction and main results}\label{S:main}

The main object of this paper is the magnetic Schr\"odinger operator
with the Dirichlet  boundary conditions in an open set $\Om\subset\R^n$.
This operator has the form

\begin{equation}\label{E:Schroed}
H_{a,V}=-\sum_{j=1}^n \left(\frac{\pa}{\pa {x^j}}+ia_j\right)^2+V,
\end{equation}
where $a_j=a_j(x)$, $V=V(x)$, $x=(x^1,\dots,x^n)\in\Om$.
We assume that $a_j$ and $V$ are real-valued functions,
$V\ge 0$,  $V\in L^1_{loc}(\Om)$ and $a\in L^\infty_{loc}(\bOm)$, i.e. for every  $j=1,\dots,n$
we have $a_j\in L^\infty_{loc}(\Om)$ and the extension of $a_j$ by $0$
to $\R^n\setminus\Om$ is in  $L^\infty_{loc}(\R^n)$.

Denote also
\begin{equation*}\label{E:nabla_a}
\nabla_a u=\nabla u +iau=
\left(\frac{\pa u}{\pa x^1}+ia_1 u,\dots, \frac{\pa u}{\pa x^n}+ia_n u\right),
\end{equation*}
and   define the quadratic form 
\begin{equation}\label{E:haV}
h_{a,V}(u,u)=\int_{\Om}(|\nabla_a u|^2+V|u|^2)dx
\end{equation}
on functions $u\in C_0^\infty(\Om)$. 
Then we can define the operator $H_{a,V}$ 
by the closure of this quadratic form in $L^2(\Om)$. 
This closure exists \cite{Leinfelder-Simader}.

In this paper we will discuss criteria for the discreteness of spectrum
and strict possitivily of $H_{a,V}$ in $L^2(\Om)$, as well as two-sided estimates
for the bottoms of the spectrum and essential spectrum. For the discreteness of spectrum
criteria the requirement $V\ge 0$ can be replaced by the semi-boundedness of $V$ from below.

We will say that  $H_{a,V}$
has a \textit{discrete spectrum} if 
its spectrum consists of isolated eigenvalues of finite
multiplicities.  It follows that the only accumulation point of
these eigenvalues is $+\infty$. Equivalently we may say that $H_{a,V}$  
has a compact resolvent.

Our first goal is to provide a direct necessary and sufficient condition  for the discreteness of the
spectrum of $H_{a,V}$. Another condition was established in \cite{Kondratiev-Mazya-Shubin}.
The difference between them is that the condition in \cite{Kondratiev-Mazya-Shubin} almost completely
separates electric and magnetic fields, whereas the condition in this paper, which is otherwise simpler,
is based on a combined characteristic of the fields. Both conditions use Wiener capacity
as was originated by A.~Molchanov \cite{Molchanov} and later developed by V.~Maz'ya
(see \cite{Mazya8}) for the usual Schr\"odinger operator without magnetic field.
(See also \cite{Mazya-Shubin, Mazya-Shubin2} for more recent results and references.)

We will denote the Wiener capacity of a compact set $F\subset \R^n$ by $\capa(F)$.
By $Q_d$ we will denote a closed cube with the edges of length $d$, which 
are parallel to coordinate axes. In case $n=2$ the capacity of $F\subset Q_d$ will be taken 
with respect to $\overset{\circ}Q_{2d}$, the interior of $Q_{2d}$.

Let $P\in \C[x]$, i.e. $P$ is a polynomial in $x^1,\dots,x^n$ with complex coefficients.
We will call it \textit{generic} on $Q_d$ if $0$ is not a critical value of the map 
$P:{\cal U}\to\C\cong \R^2$ for a neighborhood ${\cal U}$ of $Q_d$ in $\R^n$.
This means that the gradients of ${\rm Re}\, P$ and ${\rm Im}\, P$ are linearly independent
on the null-set of $P$, i.e on the set
\begin{equation*}
{\cal U}\cap P^{-1}(0)=\{x|\, x\in {\cal U},\, P(x)=0\}.
\end{equation*} 
Then this set is a non-singular algebraic submanifold of real codimension 2.
It follows that its  capacity is $0$ (see e.g. \cite{Meyers, Mazya-Havin}).
By the Sard Lemma, for any given polynomial $P$, the polynomial $P+c$ will be
generic on $Q_d$ for almost all $c\in\C$. It follows that the set
of generic (on $Q_d$) polynomials of a fixed degree is an open and dense set
in the set of all polynomials of this degree.

It is clear from the remarks above that generic polynomials form a dense set
in $C^\infty(Q_d)$.

Now we will formulate our main discreteness of spectrum result.

\begin{theorem}\label{T:discr2} 
1) Let $d_0>0$, and let $\ga=\ga(d)$ be defined for
$d\in(0,d_0)$, take values in $(0,1)$
and satisfy the condition
\begin{equation}\label{E:ga-cond}
\limsup_{d\downarrow 0}d^{-2}\ga(d)=+\infty.
\end{equation}
Then  the spectrum of $H_{a,V}$ in $L^2(\Om)$ is discrete  if and only if
for every $d\in (0,d_0)$ 
\begin{equation}\label{E:discr2} 
\inf_{F,\om}\int_{Q_d\setminus F}\left(\left|\frac{\nabla\om(x)}{i\om(x)}+a(x)\right|^2+V(x)\right)dx \to
+\infty\  {\rm as}\ Q_d\to\infty,
\end{equation}
where $F$ runs over compact sets, such that 
\begin{equation}\label{E:in-F-out}
Q_d\setminus\Om\subset F\subset Q_d
\end{equation} 
and $F$ 
satisfies the  negligibility condition
\begin{equation}\label{E:negligible}
\capa(F)\le\ga(d)\capa(Q_d),
\end{equation}
$\om$ runs over
complex valued functions from
$C^\infty(Q_d\setminus e)$,  where $e$ is a compact set 
such that $e\subset \Int F$
($e$ depends on $\om$), $|\om(x)|=1$ for every $x\in Q_d\setminus e$. Here $Q_d\to\infty$
means that the cube $Q_d$ goes to infinity with fixed $d$.

2) The same is true if 
$\om$ runs over all functions of the form $P/|P|$ where $P$ is a generic polynomial 
with complex coefficients on $Q_d$, defined on $Q_d\setminus P^{-1}(0)$. (In this case
the test functions $\om$ and sets $F$ are taken independently.)

The conditions above taken with different functions $\ga$ are equivalent.
\end{theorem}

\begin{remark}\label{R:Qd-Om}
{\rm
If for a fixed cube $Q_d$ there are no compact sets $F$ satisfying
the conditions \eqref{E:in-F-out} and \eqref{E:negligible}, then the infimum 
in \eqref{E:discr2} is naturally declared to be $+\infty$. In other words,
in \eqref{E:discr2} we can restrict ourselves to cubes $Q_d$
which are essentially in $\Om$ i.e. have negligible intersection with $\R^n\setminus \Om$
in the sense that 
\begin{equation*}
\capa(Q_d\setminus\Om)\le \ga(d)\capa(Q_d).
\end{equation*}
}
\end{remark}

\begin{remark}\label{R:infinities} 
{\rm 
The advantage of the second part of Theorem \ref{T:discr2}, as compared with the first one, 
is that the test sets $F$  and functions
$\om$ run over the corresponding families independently. However it has a disadvantage too:
the integral in \eqref{E:discr2} can be $+\infty$ if 
$P^{-1}(0)\cap(Q_d\setminus \Int(F))\ne \emptyset$.
It is easy to see that it is indeed $+\infty$ if 
$P^{-1}(0)\cap\Int(Q_d\setminus F)\ne \emptyset$. 
So this situation can be excluded from the consideration in \eqref{E:discr2}
since we are only interested in the infimum. However the integral can also be $+\infty$
if $P^{-1}(0)\cap\Int(Q_d\setminus F)=\emptyset$, but $P^{-1}(0)$ intersects with 
the boundary of $Q_d\setminus F$.
}
\end{remark}

\begin{remark}\label{R:gauge}
{\rm
We can locally represent $\om$ in the form $\om=e^{i\phi}$, where $\phi$ 
is a locally defined real-valued $C^\infty$ function. We can in fact consider $\phi$ 
as a globally  defined function with values in $\R/2\pi\Z$. Then $\nabla\phi=\nabla\om/(i\om)$,
and we see that taking infimum over $\om$'s in \eqref{E:discr2} is a way of minimizing
over different gauges. 
} 
\end{remark}

Theorem \ref{T:discr2} extends the main result of \cite{Mazya-Shubin} to the case 
of magnetic Schr\"odinger operators. (To get the discreteness of spectrum result 
of \cite{Mazya-Shubin}  we can take $a\equiv 0$ and observe that then $\om\equiv 1$  
minimizes the integral in the left hand side of \eqref{E:discr2}.)  
More general test bodies were considered in
\cite{Mazya-Shubin}  instead of cubes, and $V$ was allowed to be a positive Radon measure,
which is absolutely continuous with respect to the Wiener capacity 
(instead of a positive locally integrable function). These
results can be also extended to the magnetic Schr\"odinger  operators  
without additional difficulties. We have chosen the
simplest case for our exposition  to make it more transparent.

Another necessary an sufficient condition of the discreteness of spectrum
for the magnetic Schr\"odinger operators, obtained in \cite{Kondratiev-Mazya-Shubin},
is different in the nature of characterization of magnetic fields, and it allows
only small negligible sets (i.e. small $\ga$'s), so it does not extend the result
of \cite{Mazya-Shubin}, unlike Theorem \ref{T:discr2}. The same applies to the positivity 
results below.

Now we will formulate our main positivity result. We will say that the operator $H_{a,V}$
is {\it strictly positive} if its spectrum is contained in $[\la,+\infty)$ for some $\la>0$.
This is equivalent to the estimate
\begin{equation}\label{E:la-pos}
\int_{\Om}|u|^2 dx \le \la^{-1}h_{a,V}(u,u), \quad u\in C_0^\infty(\Om). 
\end{equation} 

\begin{theorem}\label{T:pos} Let us choose an arbitrary $\ga\in (0,1)$.
The operator $H_{a,V}$ with $V\ge 0$ is strictly positive if and only if
the following condition is satisfied: there exist $d>0$ and $\ka>0$ such that for every $Q_d$
\begin{equation}\label{E:pos}
d^{-n}\inf_{F,\om}\int_{Q_d\setminus F}\left(\left|\frac{\nabla\om(x)}{i\om(x)}+a(x)\right|^2+V(x)\right)dx\ge\ka,
\end{equation}
where $F$ runs over compact sets, such that $Q_d\setminus\Om\subset F\subset Q_d$ and $F$ 
satisfies the  negligibility condition
$\capa(F)\le\ga\capa(Q_d)$,
$\om$ runs over functions from
$C^\infty(Q_d\setminus e)$,  $e\Subset \Int F$, $|\om(x)|=1$ for all $x\in Q_d\setminus e$.

Conditions on $\om$ can be replaced by saying that $\om=P/|P|$ where $P$ 
is a generic polynomial in $Q_d$ (and then $\om$ is taken independently of 
$F$).

If $H_{a,V}$ is strictly positive, then in both cases the condition \eqref{E:pos}
is in fact satisfied for all sufficiently large $d$ (with the same $\ka$). 
\end{theorem} 

As in Theorem \ref{T:discr2}, we declare the infimum in \eqref{E:pos} to be $+\infty$
if there are no $F$'s satisfying the conditions above. 

In Sections \ref{S:two-sided} and \ref{S:2-sided-ess} we will establish two-sided estimates for 
the bottoms of the spectrum and the essential spectrum of $H_{a,V}$ in terms of 
capacitary interior diameter. These results imply Theorem \ref{T:pos}
and a weaker version of Theorem \ref{T:discr2}, with $\ga$ independent of $d$.
They extend and improve earlier results by V.~Maz'ya \cite{Mazya-74}, 
V.~Maz'ya and M.~Otelbaev \cite{Mazya-Otelbaev} 
(see also Sect. 12.2 and 12.3 in \cite{Mazya8}).

For the Dirichlet Laplacian $H_{0,0}=-\De$ in domains $\Om\subset \R^n$, stronger 
two-sided estimates  for the bottom of spectrum and essential spectrum
(with explicit constants in the estimates) were obtained in \cite{Mazya-Shubin2}
(see also more references and history there). These estimates are given
in terms of the interior capacitary radius.  

In the last Section \ref{S:example} we provide a special class
of magnetic Schr\"odinger operators where 
positivity and two-sided estimates for the bottom of the spectrum
can be found in more explicit terms: they reduce to a direct integral
of the Schr\"odinger operators without magnetic field.

\section{Sufficiency}\label{S:suff}

In this section we will prove the sufficiency of the conditions, formulated
in Theorem \ref{T:discr2}  for the discreteness of spectrum.
We will start with some general preliminaries.

For any $u\in C^\infty(Q_d)$ denote $\om=u/|u|$, which is a complex-valued function defined 
on the open set $Q_d\setminus u^{-1}(0)=\{x\in Q_d|\; u(x)\ne 0\}$, so $|\om|=1$ on this set. 

Since $\nabla_a u= \om (\nabla |u|) + |u| \nabla \om+i|u|\om a$, we have 
on $Q_d\setminus u^{-1}(0)$
\begin{equation*}
|\nabla_a u|^2=\left|\nabla |u|+i|u|\left(\frac{\nabla\om}{i\om}+a\right)\right|^2
=|\nabla|u||^2+\left|\left(\frac{\nabla\om}{i\om}+a\right)\right|^2|u|^2
\end{equation*}
because $\nabla\om/(i\om)$ is real. Therefore,
\begin{align}\label{E:haV-2}
&h_{a,V}(u,u)_{Q_d}=\int_{Q_d}(|\nabla_au|^2+V|u|^2)dx=
\int_{Q_d\setminus u^{-1}(0)}(|\nabla_au|^2+V|u|^2)dx\\
&=\int_{Q_d\setminus u^{-1}(0)}\left(|\nabla|u||^2+\left(\left|\left(\frac{\nabla\om}{i\om}+
a\right)\right|^2+V\right)|u|^2\right)dx, \notag
\end{align}
where we took into account that $\nabla u=\nabla |u|=\nabla_a u=0$ almost everywhere on $u^{-1}(0)$.
We see that  the function 
\begin{equation}\label{E:effective}
\tV=\tV[\om;a,V]=\left|\left(\frac{\nabla\om}{i\om}+a\right)\right|^2+V
\end{equation}
plays the role of an ``effective potential". 
It  is defined on $Q_d\setminus u^{-1}(0)$.

Now we will apply arguments from  Chapter 12 in \cite{Mazya8} and also Section 3 in
\cite{Kondratiev-Mazya-Shubin}.
Standard compactness arguments show that to prove the discreteness of spectrum of
$H_{a,V}$ it suffices to establish that for every $\eps>0$ there exist $d=d(\eps)$, $R=R(\eps)$
such that for every cube $Q_d$ with $\dist(Q_d,0)\ge R$ 
\begin{equation}\label{E:norm-haV}
\|u\|^2_{L^2(Q_d)}\le \eps h_{a,V}(u,u)_{Q_d}, \quad u\in C_0^\infty(\Om).
\end{equation}

The requirement $u\in C_0^\infty(\Om)$ can be replaced by $u\in \Lip_c(\Om)$
(the set of Lipschitz functions with a compact support in $\Om$) or by $u\in H_c^1(\Om)$
(the Sobolev space of functions $u\in L^2(\Om)$ with $\nabla u\in L^2(\Om)$ and with 
a compact support in $\Om$). Actually we do not need $u$ to be defined in $\Om$:
in each of the cases above we can equivalently consider $u$'s which are defined on $Q_d$ 
and vanish in a neighborhood of $Q_d\setminus\Om$. 

Below we will need the following lemma, which was proved by Maz'ya \cite{Mazya-63}
(even in a more general case of a higher order analogue of the Dirichlet integral)
and is a particular case of much more general Theorem 10.1.2, part 1, 
in \cite{Mazya8} (see also Lemma 2.1 in \cite{Kondratiev-Mazya-Shubin} 
or Lemma 4.1 in \cite{Mazya-Shubin} for simplified expositions, 
and Lemma 3.1 in \cite{Mazya-Shubin2}
for a version with an explicit constant).
 
\begin{lemma}\label{L:Mazya1}
There exists $C_n>0$ such that the following inequality holds for every complex-valued
function $u\in \Lip(Q_d)$  which
vanishes on a compact set $F\subset Q_d$ (but is not identically $0$ on $Q_d$):
\begin{equation}\label{E:cap-above}
\capa(F)\le \frac{C_n\int_{Q_d} |\nabla u(x)|^2dx}{d^{-n}\int_{Q_d} |u(x)|^2 dx}\;.
\end{equation}
\end{lemma}

The following Lemma is easily extracted from  the arguments given 
in the proof of Lemma 12.1.1  in \cite{Mazya8}.
\begin{lemma}\label{L:Mazya2}
 Let $v\in\Lip(Q_d)$, $\tilde V\in L^1_{loc}(Q_d\setminus v^{-1}(0))$, 
$\tilde V\ge 0$. Then
\begin{equation}\label{E:Mazya(2)}
\int_{Q_d}|v|^2 dx\le \min\left\{\frac{4C_nd^n}{\capa(Q_d\setminus \cM_\tau)}\int_{Q_d}|\nabla v|^2 dx,
\frac{4d^n}{\int_{\cM_\tau}\tilde V dx}\
\int_{\cM_\tau}\tilde V|v|^2 dx\right\},
\end{equation}
where $\tau=\left(\frac{1}{4d^n}\int_{Q_d}|v|^2 dx\right)^{1/2}$, $\cM_\tau=\{x|\;|v(x)|>\tau\}$,
and $C_n$ is the constant from \eqref{E:cap-above}.
(The last term in \eqref{E:Mazya(2)} is declared to be $+\infty$ if its denominator vanishes.)
\end{lemma}

\ms
{\bf Proof of Lemma \ref{L:Mazya2}.} 
Since
\begin{equation*}
|v|^2\le 2\tau^2+2(|v|-\tau)^2 \quad \hbox{on}\ {\cal M}_\tau,
\end{equation*}
we have 
\begin{equation*}
\int_{Q_d} |v|^2dx\le 2\tau^2d^n+
2\int_{\cal M_\tau} (|v|-\tau)^2dx.
\end{equation*}
Therefore
\begin{equation}\label{E:(3)}
\int_{Q_d} |v|^2 dx\le
4\int_{\cal M_\tau} (|v|-\tau)^2dx. 
\end{equation}
Using \eqref{E:(3)} and applying Lemma \ref{L:Mazya1} to the function $u=(|v|-\tau)_+$, which
equals $|v|-\tau$ on ${\cal M}_\tau$ and $0$ on
$Q_d\setminus {\cal M}_\tau$, we see that
\begin{equation*}
\capa(Q_d\setminus{\cal M}_\tau)\le 
\frac{C_n\int_{{\cal M}_\tau}|\nabla(|v|-\tau)|^2 dx}
{d^{-n}\int_{Q_d}|u|^2 dx}\le
\frac{4C_n\int_{Q_d}|\nabla v|^2 dx}
{d^{-n}\int_{Q_d}|v|^2 dx}\;.
\end{equation*}
where $C_n$ is the constant from \eqref{E:cap-above}. 
Therefore
\begin{equation}\label{E:(4)}
\int_{Q_d}|v|^2 dx \le\frac{4C_n d^n\int_{Q_d}|\nabla v|^2 dx}{\capa(Q_d\setminus{\cal M}_\tau)}\;.
\end{equation}
On the other hand,
\begin{equation*}
\int_{{\cal M}_\tau} \tV |v|^2  dx\ge
\tau^2\int_{{\cal M}_\tau} \tV dx
=\frac{1}{4d^n}\int_{Q_d}|v|^2dx\cdot \int_{{\cal M}_\tau} \tV dx,
\end{equation*}
hence
\begin{equation}\label{E:(5)}
\int_{Q_d} |v|^2 dx\le 
\frac{4 d^n}{ \int_{{\cal M}_\tau} \tV dx}\int_{\cal M_\tau} \tV|v|^2 dx,
\end{equation}
where the right hand side is declared to be $+\infty$ if the denominator is $0$.

The resulting inequality \eqref{E:Mazya(2)} follows from \eqref{E:(4)} and \eqref{E:(5)}. $\square$

\begin{corollary}\label{C:Mazya2} There exists $C_n>0$ such that the following holds. 
Assume that $v\in\Lip(Q_d)$, $\tilde V\in L^1_{loc}(Q_d\setminus v^{-1}(0))$, 
$\tilde V\ge 0$. Then for every $\ga\in(0,1)$
\begin{equation}\label{E:Mazya(2)-cor}
\int_{Q_d}|v|^2 dx\le \frac{C_nd^2}{\gamma}  \int_{Q_d}|\nabla v|^2 dx
+\frac{4d^n}{\underset{F}\inf\int_{Q_d\setminus F}\tilde V dx}\
\int_{Q_d\setminus v^{-1}(0)}\tilde V|v|^2 dx,
\end{equation}
where the infimum is taken over all compact sets $F\subset Q_d$ such that
$\Int F\supset v^{-1}(0)$ and $F$ satisfies 
the negligibility condition \eqref{E:negligible}. If there is no such $F$'s, 
then we declare the infimum to be $+\infty$ and the last term 
itself to be $0$.
The last term in \eqref{E:Mazya(2)-cor} is declared to be $+\infty$ if the infimum vanishes.
\end{corollary}

{\bf Proof.} Assume first that $\capa(Q_d\setminus \cM_\tau)>\ga\capa(Q_d)$. 
Then Lemma \ref{L:Mazya2} implies that 
\begin{equation}\label{E:part1}
\int_{Q_d}|v|^2 dx\le \frac{C_n^{(1)}d^n}{\ga\capa(Q_d)}\int_{Q_d}|\nabla v|^2 dx=
\frac{C_n^{(2)}d^2}{\ga}\int_{Q_d}|\nabla v|^2 dx. 
\end{equation}
In the opposite case $\capa(Q_d\setminus \cM_\tau)\le \ga\capa(Q_d)$ we
can use the second term in the braces in the right-hand side of \eqref{E:Mazya(2)} and
replace the integral in the denominator by the infimum of such integrals over $Q_d\setminus F$
with the conditions on $F$ as formulated in the Corollary. (Note that in this case
the family of admissible $F$'s is not empty because it includes $Q_d\setminus\cM_\tau$.) 
We get then
\begin{equation}\label{E:part2}
\int_{Q_d}|v|^2 dx\le \frac{4d^n}{\underset{F}\inf\int_{Q_d\setminus F}\tilde V dx}\
\int_{Q_d\setminus v^{-1}(0)}\tilde V|v|^2 dx.
\end{equation}
Combining \eqref{E:part1} and \eqref{E:part2}, we get \eqref{E:Mazya(2)-cor} 
for arbitrary $\ga\in (0,1)$. $\square$ 

\ms
{\bf Proof of sufficiency in Theorem \ref{T:discr2}.}
1) Let us apply \eqref{E:Mazya(2)-cor} with $v=|u|$ where $u\in C_0^\infty(\Om)$ 
and $\tilde V$ given by
\eqref{E:effective} with $\om=u/|u|$. Note that the condition $\Int F\supset v^{-1}(0)$
implies that $\Int F\supset Q_d\setminus \Om$ because $\supp v$ is a compact subset in $\Om$.
We see  that to achieve
\eqref{E:norm-haV}  it suffices that the following two conditions
are satisfied (for the same cube $Q_d$)
\begin{equation}\label{E:2-conditions}
\frac{C_n d^2}{\ga}\le \eps,\qquad \frac{4d^n}{\underset{F}\inf\int_{Q_d\setminus F}\tilde V dx}\le\eps,
\end{equation}
where $\gamma$ in the first inequality and in the negligibility condition \eqref{E:negligible}
may depend upon $Q_d$ (and even  upon $a$ and $V$ as well). Clearly, these inequalities
will hold if we require that the first inequality holds together with the second one replaced
by the stronger inequality
\begin{equation}\label{E:stronger-ineq}
\frac{4d^n}{\underset{F,\om}\inf\int_{Q_d\setminus F}\tilde V dx}\le\eps.
\end{equation}
Now we can take $\ga=\ga(d)$ satisfying \eqref{E:ga-cond} to see that it is sufficient that for for every
$\eps>0$ there exists $d\in(0,d_0)$ and $R=R(d)>0$ such that for every cube $Q_d$ with $\dist(Q_d,0)\ge R$
\begin{equation}\label{E:d-conditions}
C_n d^2 \ga(d)^{-1}\le\eps, \qquad 
\underset{F,\om}\inf\;\frac{1}{d^n}\int_{Q_d\setminus F}\tilde V dx \ge 4\eps^{-1}. 
\end{equation}
It is clear from \eqref{E:ga-cond} that   
we can choose $\eps=\max\{C_n,4\}d^2\ga(d)^{-1}$ along a sequence $d=d_k\to 0$ to conclude that
the conditions \eqref{E:d-conditions} can be replaced by a single condition
\begin{equation}\label{E:discr3}
\underset{F,\om}\inf\;\frac{1}{d^n}\int_{Q_d\setminus F}\tilde Vdx \ge d^{-2}\ga(d), \quad \tV=\tV[\om;a,V],
\end{equation}
which should be satisfied for all $d\in (0,d_0)$ and for distant cubes,
i.e. when $\dist(Q_d,0)\ge R$, $R=R(d)$. This proves the sufficiency in the first part of
Theorem \ref{T:discr2}, because \eqref{E:discr2} obviously implies \eqref{E:discr3}. 

\medskip
2) Let us prove sufficiency in the second part of Theorem \ref{T:discr2}.
We will argue as in the first
part of this proof. 
Since the generic polynomials are dense in $C^{\infty}(Q_d)$ 
(in its standard Frechet topology), 
it suffices to establish that for every $\eps>0$ the estimate \eqref{E:norm-haV} holds
for distant cubes $Q_d$. It follows that it is sufficient to require that
there exists $d_0>0$ such that
for every $d\in (0, d_0)$ 
the inequality  \eqref{E:discr3} holds for distant cubes $Q_d$ 
(i.e. cubes with $\dist(Q_d,0)\ge R(d)$),
for generic polynomials $u$ on $Q_d$,
and with $\tV=\tV[\om;a,V]$, where $\om=P/|P|$, $P$ is another
generic polynomial on $Q_d$, such that $P^{-1}(0)\cap Q_d$ is in the interior of $F$ 
with respect to $Q_d$, i.e. 
\begin{equation}\label{E:dist}
\dist(P^{-1}(0)\cap Q_d, Q_d\setminus F)>0.
\end{equation}
According to the formulation of Theorem \ref{T:discr2}, part 2, we can  assume that \eqref{E:discr3} 
is fulfilled if
we  allow arbitrary $\om=P/|P|$, unrelated to $F$, i.e. if we abandon the condition \eqref{E:dist}.
But then the infimum in \eqref{E:discr3} becomes smaller and we get a stronger condition. So both versions
of \eqref{E:discr3} are satisfied, which ends the proof of sufficiency in Theorem \ref{T:discr2}.
$\square$

\begin{remark}\label{R:inf-equiv} 
{\rm
Let us denote the infimum in \eqref{E:discr3} by $I(\ga)$, where we assume that $\om=P/|P|$ and $F$
(satisfying \eqref{E:in-F-out} and \eqref{E:negligible})
are taken independently of each other. By 
$\tI(\ga)$ we will denote the same infimum but with additional condition \eqref{E:dist} imposed upon 
$P$  and $F$. Clearly, $I(\ga)\le\tI(\ga)$. 

On the other hand  
the inequality $\tI(\ga)\le I(\ga')$ holds for any positive $\ga'<\ga$. This would follow if
we prove that for every fixed $\om=P/|P|$, where $P$ is a generic polynomial on $Q_d$, and for every
compact $F'$, $Q_d\setminus\Om\subset F'\subset Q_d$, with $\capa(F')\le\ga'\capa(Q_d)$, 
there exists a compact $F$, $Q_d\setminus\Om\subset F\subset Q_d$, with $\capa(F)\le\ga\capa(Q_d)$, 
satisfying \eqref{E:dist}, such that
\begin{equation*}
\int_{Q_d\setminus F}\tV[\om;a,V] dx \le \int_{Q_d\setminus F'}\tV[\om;a,V] dx,
\end{equation*}
with the same $\om=P/|P|$. The last inequality will be fulfilled automatically if $F\supset F'$.
So we can take $F=F'\cup U$ where $U$ is the closure of a sufficiently small neighborhood of 
$P^{-1}(0)\cap Q_d$ in $Q_d$. Then the condition \eqref{E:dist} will be  satisfied. Besides,
we can choose $U$ to have an arbitrarily small capacity because $\capa(P^{-1}(0))=0$. Then the inequality
$\capa(F)\le\ga\capa(Q_d)$ will also hold due to subadditivity of capacity. 

This argument shows that the condition \eqref{E:discr3} with independent $F$ and $\om=P/|P|$
follows from the same condition with the additional requirement \eqref{E:dist}, but with an arbitrary
smaller $\ga$. (We can take e.g. $(1-\ka)\ga(d)$ instead of $\ga(d)$, where $\ka>0$ is arbitrarily small.)
So these conditions are almost equivalent. 
In particular, the corresponding conditions \eqref{E:discr2} are equivalent. 
}
\end{remark}

\begin{remark}
{\rm
Instead of the condition \eqref{E:discr2}, it is sufficient to require
a weaker condition \eqref{E:discr3}. (The left hand side of \eqref{E:discr3} is not
required to tend to $+\infty$ as $Q_d\to\infty$.)  But, 
as we will see later, the stronger condition \eqref{E:discr2} is also necessary, so
the conditions \eqref{E:discr2} and \eqref{E:discr3} are in fact equivalent
for any $\gamma$, satisfying  \eqref{E:ga-cond}.
In particular, it would follow that the above versions of the condition \eqref{E:discr3} are equaivalent.
}
\end{remark}

\section{Necessity}\label{S:nec}

In this section we will prove the necessity of the conditions formulated 
in Theorem \ref{T:discr2} for the discreteness of spectrum.

\ms
1) Let us assume that $H_{a,V}$ has a discrete spectrum. This implies that for every $d>0$
\begin{equation}\label{E:nec-loc}
\inf_u\;\frac{h_{a,V}(u,u)_{Q_d}}{\|u\|^2_{L^2(Q_d)}}\to +\infty \quad {\rm as} \quad Q_d\to\infty,
\end{equation}
where the infimum is taken over all $u\in\Lip(Q_d)$, $u\not\equiv 0$,
$u=0$ in a neighborhood of $Q_d\setminus \Om$
(see e.g. Theorem 1.2 in \cite{Kondratiev-Shubin}). Let us fix $d>0$. 
Then \eqref{E:nec-loc} means that  for every 
$\eps>0$ there exists $R=R(\eps)$ that for every cube $Q_d$ with ${\rm dist}(Q_d,0)\ge R$
\begin{equation}\label{E:nec-loc2}
\|u\|^2_{L^2(Q_d)}\le \eps h_{a,V}(u,u)_{Q_d}, \quad u\in \Lip(Q_d),\ 
(\supp u)\cap (Q_d\setminus\Om)=\emptyset. 
\end{equation}
We would like to use this estimate with  a test function $u=(1-P_F)\om$. 
Here $F\subset\R^n$, $\Int F\supset Q_d\setminus\Om$, $F$ is a {\it regular} compact set, 
i.e. $F$ is a compact subset in the cube
$Q_{3d/2}$ (with the same center as $Q_d$) and $F$ is the closure of an open set with a smooth boundary;
$P_F$ is the  equilibrium potential of $F$, i.e. for $n\ge 3$ we have 
$P_F\in C(\R^n)$, $P_F=1$ on $F$, $\De P_F=0$ on $\R^n\setminus F$, $P_F(x)\to 0$ as $|x|\to\infty$,
and for $n=2$ we have $P_F\in C(\R^2)$, $P_F=1$ on $F$, $\De P_F=0$ on $(\Int Q_{2d})\setminus F$,
$P_F=0$ on $\R^2\setminus Q_{2d}$;  
$\om\in C^\infty(Q_d\setminus e)$,
$|\om(x)|=1$ for all $x\in Q_d\setminus e$, $e\subset \Int F$ and $F,\om$ are chosen so that 
\begin{equation}\label{E:near-inf}
\int_{Q_d\setminus F} \tV dx \le \inf_{F,\om} \int_{Q_d\setminus F} \tV dx + \de d^{n-2},
\end{equation}
where we use the same notations as in Sect. \ref{S:suff} (in particular
$\tV=\tV[\om;a,V]$ is given by \eqref{E:effective}), and $\de>0$ is sufficiently small. 
It is well known that $0\le P_F\le 1$ everywhere and
\begin{equation}\label{E:phi-F}
\int_{\R^n}|\nabla P_F|^2 dx=\capa(F).
\end{equation}

Clearly, $|u|=1-P_F$ and $\om(x)=u(x)/|u(x)|$ if $u(x)\ne 0$. Therefore $u(x)=|u(x)|\om(x)$
for all $x\in Q_d\setminus e$. The calculations in Sect. \ref{S:suff} are applicable in this case
and lead to the formula \eqref{E:haV-2} and to the ``effective" potential $\tV$ of the same form
\eqref{E:effective}, defined on $Q_d\setminus e$.  It follows that
\begin{align}\label{E:haV-above}
&h_{a,V}(u,u)_{Q_d}=
\int_{Q_d}|\nabla P_F|^2 dx
+\int_{Q_d}\left(\left|\frac{\nabla\om}{i\om}+a\right|^2+V\right)|1-P_F|^2 dx\\
&\le  \capa(F)  + \int_{Q_d\setminus F}\left(\left|\frac{\nabla\om}{i\om}+a\right|^2+V\right) dx
= \capa(F)+ \int_{Q_d\setminus F}\tV dx \notag\\
& \le \capa(F) + \de d^{n-2} + \inf_{F,\om}\int_{Q_d\setminus F}\tV dx. \notag
\end{align}

Now we need to estimate  $\|u\|_{L^2(Q_d)}$ from below.
To this end we need the following
\begin{lemma}\label{L:L2-norm-below}
There exists $C=C_n>0$ such that for every $\eta\in(0,1/2]$
\begin{equation}\label{E:L2-norm-below}
\left(1-\frac{\capa(F)}{\capa(Q_d)}\right)^2
\le C_n\left[\eta\frac{\capa(F)}{\capa(Q_d)}+\eta^{-1} d^{-n}\int_{Q_d}(1-P_F)^2 dx\right].
\end{equation}
\end{lemma}

For the proof of this Lemma see  \cite{Mazya-Shubin} (formula (3.10) there).  

\medskip
Assuming that $\capa(F)\le\ga\capa(Q_d)$ with $\ga\in(0,1)$, we obtain 
\begin{equation}\label{E:L2-norm-below2}
(1-\ga)^2\le C_n\left[\eta\ga+\eta^{-1} d^{-n}\int_{Q_d}(1-P_F)^2 dx\right].
\end{equation}
Choosing
\begin{equation}\label{E:eta}
\eta=\min\left\{\frac{1}{2}, \frac{(1-\ga)^2}{2\ga C_n}\right\},
\end{equation} 
we obtain 
\begin{equation}\label{E:int-below}
\int_{Q_d}(1-P_F)^2 dx\ge (2C_n)^{-1}\eta(1-\ga)^2 d^n.
\end{equation}

Taking into account that $|u|=1-P_F$ and 
using the estimates  
\eqref{E:haV-above} and \eqref{E:int-below} in \eqref{E:nec-loc2}, we see
that for distant cubes
\begin{equation}\label{E:final-nec}
(2C_n)^{-1}\eta(1-\ga)^2 d^n\le \eps \left(\capa(F) + \de d^{n-2} 
+ \inf_{F,\om}\int_{Q_d\setminus F}\tV dx\right),
\end{equation}
where the infimum is taken over all regular $F$ satisfying $\capa(F)\le\ga\capa(Q_d)$. 
But now we can approximate an arbitrary $F$, satisfying the same inequality, by 
regular sets from above, using the well-known continuity property of the capacity
(see e.g. \cite{Mazya8}, Sect. 2.2.1). Then we obtain the same inequality with the infimum taken
over arbitrary (not necessarily regular) negligible compact sets $F$. 
It follows that \eqref{E:discr2} is satisfied, which ends 
proof of the first part of Theorem \ref{T:discr2}. 

\ms
2) According to the first part of Theorem \ref{T:discr2},
which we already established, the discreteness of spectrum implies that the condition \eqref{E:discr2}
is fulfilled if the infimum taken over $\om, F$ such that $\om=P/|P|$ where $P$ is a generic 
polynomial on $Q_d$, $F$ is a compact subset in $Q_d$ and the condition \eqref{E:dist}
is satisfied. It remains to get rid of the condition \eqref{E:dist}. This is easily done
by the same arguments as in the proof of sufficiency in Theorem \ref{T:discr2} 
and in Remark \ref{R:inf-equiv}.
This ends the proof of Theorem \ref{T:discr2}.  $\square$

\section{Positivity}\label{S:pos}

In this section we will prove  Theorem \ref{T:pos} and provide its interesting corollary.

\ms
\textbf{Proof of sufficiency in Theorem \ref{T:pos}.}  Let us assume that \eqref{E:pos} holds for some
$d>0$ and $\ka>0$ (with any of two versions for the set of $\om$'s). 
Using Corollary \ref{C:Mazya2}, as in the proof of 
the sufficiency in Theorem \ref{T:discr2}, we come to 
the conclusion that the estimate \eqref{E:la-pos} holds with
\begin{equation}\label{E:la}
\la=\min\left\{C_n^{-1}d^{-2}\ga, \ka/4\right\},
\end{equation}
hence $H_{a,V}$ is strictly positive. $\square$

\ms
\textbf{Proof of necessity in Theorem \ref{T:pos}.} Let us assume that
$H_{a,V}$ is strictly positive, i.e. the estimate \eqref{E:la-pos} holds.
Then arguing as in the proof of necessity in Theorem \ref{T:discr2},
we come to the estimate \eqref{E:final-nec} with $\eps=\la^{-1}$ and $\eta$
given by \eqref{E:eta}. It follows that
\begin{align}\label{E:final-pos}
&d^{-n}\inf_{F,\om}\int_{Q_d\setminus F}
\left(\left|\frac{\nabla\om}{i\om}+a\right|^2+V\right)dx\\
&\ge (2C_n)^{-1}\eta(1-\ga)^{2}\la-d^{-2}(\capa(Q_1)+\de).
\notag
\end{align}
This implies  the inequality
\begin{equation}\label{E:final-pos2}
d^{-n}\inf_{F,\om}\int_{Q_d\setminus F}
\left(\left|\frac{\nabla\om}{i\om}+a\right|^2+V\right)dx
\ge (4C_n)^{-1}\eta(1-\ga)^{2}\la,
\end{equation}
provided $d>0$ is chosen so that
\begin{equation}\label{E:d-ge}
d^2\ge 4C_n\eta^{-1}(1-\ga)^{-2}\la^{-1}(\capa(Q_1)+\de).
\end{equation}
This works for both versions of the choices of $\om$'s and so ends the proof 
of Theorem \ref{T:pos}. $\square$

\section{Two-sided estimates for the bottom of the spectrum}\label{S:two-sided} 
In this section we will establish two-sided estimates for the bottom
of the spectrum for the Schr\"odinger operators
$H_{a,V}$ in $L^2(\Om)$. These estimates extend and improve results 
by V.~Maz'ya and  M.~Otelbaev
\cite{Mazya-Otelbaev} (see also \cite{Mazya8}, Sect. 12.2, 12.3)
where the case of Schr\"odinger operators without magnetic fields was considered. 
The results are based on
the notion of the capacitary interior diameter which is defined as follows: 
\begin{equation}\label{E:D}
D=D(\Om,\ga,a,V)=\sup_{Q_d}\left\{d:\; d^{n-2}\ge \inf_{F,\om}\int_{Q_d\setminus F}
\left(\left|\frac{\nabla\om}{i\om}+a\right|^2+V\right)dx\right\},
\end{equation}
where the choice of the pairs $F,\om$ is as in Theorem \ref{T:pos}
(with any of two options there), with $\ga\in (0,1)$ assumed to be a constant. 
It is easy to see  that $D>0$ (take $d$ to be very small). 
On the other hand it may happen that $D=+\infty$; for example
this is the case if $\Om$ contains arbitrarily large cubes (e.g. $\Om=\R^n$)
and both $a$ and $V$ vanish identically, i.e. when
$H_{a,V}=-\De$ in $L^2(\R^n)$. So generally $0<D\le +\infty$. 

It is easy to see that $D$ is an increasing function of $\Om$ and $\ga$
(provided $a$ and $V$ are fixed).

The definition of $D$ by \eqref{E:D} can be extended to the case when $\ga=\ga(d)$, i.e.
$\ga:(0,+\infty)\to (0,1)$, but for simplicity we will only consider the case
when $\ga$ is a constant.  

\begin{theorem}\label{T:2-sided} 
For every $\ga\in(0,1)$ there exists $C=C(\ga,n)>0$, such that
\begin{equation}\label{E:2-sided}
C^{-1}D^{-2} \le \la \le CD^{-2},
\end{equation}
where $\la$ is the bottom of the Dirichlet spectrum of $H_{a,V}$ in $L^2(\Om)$,
i.e. the best constant in \eqref{E:la-pos}.
\end{theorem}

\begin{remark} 
{\rm
Note that $H_{a,V}$ is strictly positive if and only if $\la>0$. Therefore \eqref{E:2-sided}
implies that the strict positivity is equivalent to the inequality $D<+\infty$.
More precisely, the first inequality in \eqref{E:2-sided}, estimating $\la$ from below, implies 
the sufficiency in Theorem \ref{T:pos}, whereas the second one, estimating $\la$
from above, implies the necessity in this theorem. We will obtain proofs of 
these inequalities by analyzing corresponding parts of the proof of Theorem 
\ref{T:pos}.
}
\end{remark}

{\bf Proof of Theorem \ref{T:2-sided}.} 1) Let us start with proving the first inequality, 
estimating $\la$ from below.  
Denote the bottom of Neumann spectrum of $H_{a,V}$ on $Q_d$ by $\mu(Q_d)$, i.e.
$\mu(Q_d)$ is the left hand side of \eqref{E:nec-loc}.
(Note that it depends upon $\Om$ too.)
Then \eqref{E:Mazya(2)-cor} implies that for every $d>0$ 
and  every cube $Q_d$
\begin{align*}
\mu(Q_d)^{-1}&\le C_1\max\left\{\frac{d^2}{\ga},\ \frac{d^n}{\inf_{F,\om}\int_{Q_d\setminus F}\tV dx}\right\}\\
&\le C_1\ga^{-1}\max\left\{d^{2}, \frac{d^n}{\inf_{F,\om}\int_{Q_d\setminus F}\tV dx}\right\},
\end{align*}
where $\tV=\tV[\om;a,V]$ is the ``effective potential" defined by \eqref{E:effective}.
Therefore,
\begin{equation}\label{E:mu-ge}
\mu(Q_d)\ge \ga C_1^{-1}\min\left\{d^{-2}, d^{-n}\inf_{F,\om}\int_{Q_d\setminus F}\tV dx\right\}.
\end{equation}
Let us assume that for some $d>0$ we have
\begin{equation}\label{E:d-2-le} 
d^{-2}\le d^{-n}\inf_{F,\om}\int_{Q_d\setminus F}\tV dx 
\quad\text{\rm for all}\ \  Q_d.
\end{equation}
(This holds in particular if $d>D$.)
Then we obviously have 
\begin{equation}\label{E:la-ge}
\la\ge \inf_{Q_d} \mu(Q_d)\ge \ga C_1^{-1}d^{-2} 
\end{equation}
for this particular $d$. Taking limit as $d\downarrow D$, we obtain the same inequality 
with $d=D$, which proves the left inequality in \eqref{E:2-sided}.

\ms
2) Now let us prove the second inequality in \eqref{E:2-sided}, 
estimating $\la$ from above.
To this end let us look at the inequalities \eqref{E:final-pos2} and \eqref{E:d-ge}. 
They imply that there exists $C=C(\ga)>0$ such that for every cube $Q_d$ 
at least one of the inequalities
\begin{equation}\label{E:2-ineq}
\la\le Cd^{-2}, \quad \la\le C d^{-n}\inf_{F,\om}\int_{Q_d\setminus F}\tV dx, 
\end{equation}
must hold. By definition of $D$ (see \eqref{E:D}) this implies that
$\la\le CD^{-2}$ which ends the proof.
$\square$

\begin{remark} 
{\rm For the Dirichlet Laplacian (i.e. for $H_{0,0}=-\De$) in domains $\Om$
and for small $\ga>0$ the result of Theorem \ref{T:2-sided} was proved in \cite{Mazya-74}.
This result for arbitrary $\ga\in(0,1)$ was obtained in \cite{Mazya-Shubin2}, where also 
explicit values of the constants in the upper and lower bounds were given. 
}
\end{remark}

\section{Persson type theorem for magnetic Schr\"odinger operators}\label{S:Persson}

This section contains a preparatory result, expressing the bottom of the
essential spectrum of $H_{a,V}$ in $L^2(\Om)$ as a limit of the bottoms of the Dirichlet
spectrum of this operator on the exteriors of large balls. The first result of this kind 
is probably due to Persson \cite{Persson}  (for the usual Schr\"odinger operators, without 
magnetic field), see also Chapter 3 in \cite{Agmon82} and Theorem 3.12 in \cite{CFKS}.
(In particular, the Laplacian is replaced by general second-order operators in divergence form
in \cite{Agmon82}.)
 The arguments given in these sources can be extended to our case. Nevertheless for the sake
of convenience of the reader we offer a proof which seems to be different from what we 
have seen in the literature.

For any open set $U\subset\R^n$, denote by $\la(U; H_{a,V})$ the bottom of the Dirichlet
spectrum of $H_{a,V}$ in $L^2(U)$, i.e. the spectrum of the operator defined by 
the closure of the quadratic form $h_{a,V}$ (see \eqref{E:haV}) in $L^2(U)$
from the initial domain $C_0^\infty(U)$. In other words,
\begin{equation}\label{E:la-Om}
\la(U;H_{a,V})=
\inf\left\{\left.\frac{h_{a,V}(u,u)}{(u,u)}\right|u\in C_0^{\infty}(U)\setminus\{0\}\right\},
\end{equation}
where $(\cdot,\cdot)$ means the scalar product in $L^2(U)$.
 
Usually $H_{a,V}$ will be fixed in our arguments, and in this case we will write $\la(U)$ 
instead of $\la(U; H_{a,V})$ if this does not lead to a confusion.

Note that $U\subset U'$ implies $\la(U)\ge \la(U')$.

For any self-adjoint operator $H$ in a Hilbert space $\cH$ and any $\la\in\R$ denote by $E_\la$ 
(or $E_\la(H)$) the spectral projection of $H$ corresponding  to the interval $(-\infty,\la)$.
Let us introduce the ``counting function" of the spectrum by
\begin{equation}\label{E:N-la}
N(\la)=N(\la;H)=\Tr E_\la=\dim{\rm Im} E_\la.
\end{equation}
It is an increasing function of $\la$ with values in $[0,+\infty]$. If $H$ is semibounded below,
and $\la_\infty$ is the bottom of its essential spectrum $\sigma_{ess}(H)$, then
\begin{equation}\label{E:la-infty-abstr}
\la_\infty = \sup \{\la\in\R|\; N(\la)<+\infty\}.
\end{equation} 

The following Lemma is a well known variational principle
(see e.g. \cite{Berezin-Shubin}, Section 3 in Appendix 1). In this form
it is often attributed to I.M.~Glazman.

\begin{lemma}[\rm Glazman's Lemma]\label{L:Glazman} For every $\la\in\R$,
\begin{equation}\label{E:Glazman}
N(\la)=\sup\{\dim L|\; L\subset Q(h), h(u,u)<\la(u,u), \forall u\in L\setminus\{0\}\},
\end{equation}
where $h$ is the quadratic form of $H$, $Q(h)$ is the domain of $h$, $L$ is 
a linear subspace of $Q(h)$,  $(\cdot,\cdot)$ stands for the scalar product 
in $\cH$.

This holds also if we replace $Q(h)$ by any core of $h$.
\end{lemma}

\begin{corollary}\label{C:N-la-2} For the operator $H_{a,V}$ in $L^2(\Om)$
with the Dirichlet boundary condition
\begin{equation}\label{E:N-la-2}
N(\la)=\sup\{\dim L|\; L\subset C_0^\infty(\Om), h_{a,V}(u,u)<\la(u,u), \forall u\in L\setminus\{0\}\}.
\end{equation}
\end{corollary}

The following theorem is a version of the Persson theorem \cite{Persson}
(see also Theorem 3.12 in \cite{CFKS}, as well as \cite{Agmon82}).

\begin{theorem}\label{P:la-infty} 
For $\la_\infty=\inf\sigma_{ess}(H_{a,V})$ in $L^2(\Om)$ we have
\begin{equation}\label{E:la-infty}
\la_\infty=\lim_{R\to\infty}\la(\Om\setminus \bB_R(0)),
\end{equation}
where $\bB_R(0)$ is the closed ball with the radius $R$ and center at $0$.
\end{theorem}
 
Note that the limit in \eqref{E:la-infty} exists in $[0,+\infty]$
because
$\la(\Om\setminus \bB_R(0))$ increases with respect to $R$.

\ms
\textbf{Proof of Theorem \ref{P:la-infty}}.
1) Let us prove first that $\la_\infty$ is not larger than the right hand side 
in \eqref{E:la-infty}. To this end it suffices to consider the case
when the right hand side is finite.  

Let us take an arbitrary $\tilde\lambda$ such that
\begin{equation}\label{E:tilde-la}
\tilde\lambda>\lim_{R\to\infty}\la(\Om\setminus \bB_R(0))=\sup_R \la(\Om\setminus \bB_R(0)).
\end{equation}
Taking $R_1>0$ arbitrary and using the inequality $\tilde\lambda>\la(\Om\setminus \bB_{R_1}(0))$,
we can find a function $\psi_1\in C_0^\infty(\Om\setminus \bB_{R_1}(0))$, such that $\|\psi_1\|=1$
(here $\|\cdot\|$ means the norm in $L^2(\R^n)$) and $h_{a,V}(\psi_1,\psi_1)< \tilde\lambda$.
Choose $R_2>R_1$ so that $\supp\psi_1\subset B_{R_2}(0)$. Now using the inequality 
$\tilde\lambda>\la(\Om\setminus \bB_{R_2}(0))$, we can construct 
$\psi_2\in C_0^\infty(\Om\setminus \bB_{R_2}(0))$ so that $\|\psi_2\|=1$
and $h_{a,V}(\psi_2,\psi_2)< \tilde\lambda$. Proceeding by induction, we can
construct an orthonormal sequence of functions $\psi_k\in C_0^\infty(\Om)$, $k=1,2,\dots$,
with disjoint supports, such that $h_{a,V}(\psi_k,\psi_k)<\tilde\lambda$
for all $k$. It follows that $h_{a,V}(u,u)< \tilde\lambda(u,u)$ for any $u\ne 0$ in the linear span 
$L$ of the sequence $\{\psi_k\}$. This implies that $N(\tilde\lambda)=+\infty$, 
hence $\la_\infty\le\tilde\lambda$ due to \eqref{E:la-infty-abstr}. Since $\tilde\lambda$
is an arbitrary number satisfying \eqref{E:tilde-la}, this proves that $\la_\infty$
does not exceed the right hand side of \eqref{E:la-infty}.

2) Now let us prove that $\la_\infty$ is not smaller than 
right hand side of \eqref{E:la-infty}. To this end it is sufficient to prove that
\begin{equation}\label{E:la-infty-below}
\la_\infty\ge\la(\Om\setminus K),
\end{equation}
for every compact set $K\subset\Om$. 
It is enough to consider the case when $\la_\infty<\infty$.
Ad absurdum let us assume that $\la_\infty<\la(\Om\setminus K)$
for some $K$. Let us choose $\tilde\lambda$ so that 
\begin{equation}\label{E:la-3}
\la_\infty<\tilde\lambda<\la(\Om\setminus K). 
\end{equation}
Due to \eqref{E:la-infty-abstr} we have $N(\tilde\lambda)=+\infty$,
and the same is true if we replace $\tilde\la$ by $\tilde\la-\eps$
for a sufficiently small $\eps>0$. Therefore, due to Corollary \ref{C:N-la-2}, 
for every integer $N>0$ we can find a subspace $L\subset C_0^{\infty}(\Om)$,
such that $\dim L=N$ and
\begin{equation}\label{E:haV-on-L}
h_{a,V}(u,u)\le\tilde\lambda(u,u), \quad u\in L.
\end{equation}
On the other hand, the opposite inequality is true if $\supp u\subset \Om\setminus K$.
More precisely, according to \eqref{E:la-Om},
\begin{equation}\label{E:haV-out}
h_{a,V}(u,u)\ge\lambda(\Om\setminus K)(u,u), \quad u\in C_0^\infty(\Om\setminus K).
\end{equation}
We will establish that the combination of \eqref{E:haV-on-L}  and \eqref{E:haV-out}
is impossible if $N$ is sufficiently large. To this end we will split $|u|^2$ 
for every $u\in L$
into a sum $|u_0|^2+|u_1|^2$ where $u_0$ is  supported in a neighborhood of $K$ and  $u_1$
has its support in $\Om\setminus K$. This splitting is conveniently done
by use of the IMS localization formula (see e.g. Section 3.1 in \cite{CFKS} and 
Lemma 3.1 in \cite{Shubin96}).

Let us choose $R_0>0$, so that $K\subset B_{R_0}(0)$, and take $R>2R_0$. 
We can choose functions $J_0\in C_0^\infty(B_{2R}(0))$ 
and $J_1\in C^\infty(\R^n)$, such that 
$J_0=1$ on $\bB_{R_0}(0)$, $0\le J_0(x)\le 1$ for all 
$x\in \R^n$, $J_1\ge 0$, $J_0^2+J_1^2\equiv 1$
(hence $J_1\le 1$ and $J_1=0$ on $\bB_{R_0}(0)$), and
\begin{equation}\label{E:der-est}
\sup(|\nabla J_0|+|\nabla J_1|)\le C_n R^{-1}.
\end{equation}
Due to the IMS localization formula we obtain for any $u\in C_0^\infty(\R^n)$ 
\begin{equation}\label{E:IMS}
h_{a,V}(u,u)=\sum_{k=0}^1 h_{a,V}(J_ku,J_ku)-\sum_{k=0}^1(|\nabla J_k|^2u,u). 
\end{equation}
It follows from \eqref{E:der-est} and \eqref{E:IMS} that
\begin{equation}\label{E:IMS-ge}
h_{a,V}(u,u)\ge h_{a,V}(J_0u,J_0u)+h_{a,V}(J_1u,J_1u)-C_n^2 R^{-2}(u,u).
\end{equation}
The inequality 
\eqref{E:haV-out} implies 
\begin{equation}\label{E:haV-J1}
h_{a,V}(J_1u,J_1u)\ge \la(\Om\setminus K) (J_1u,J_1u), \quad u\in C_0^\infty(\Om).
\end{equation}
Taking $u\in L$, we obtain from \eqref{E:haV-on-L}, \eqref{E:IMS-ge} and \eqref{E:haV-J1} that
\begin{align*}
&h_{a,V}(J_0u,J_0u)\le h_{a,V}(u,u)+C_n^2R^{-2}(u,u)-h_{a,V}(J_1u,J_1u)\\
&\le (\tilde\la +C_n^2R^{-2})(u,u)-\la(\Om\setminus K)(J_1u,J_1u)\\
&= (\tilde\la +C_n^2R^{-2})(J_0u,J_0u) - 
(\la(\Om\setminus K)-\tilde\la-C_n^2R^{-2})(J_1u,J_1u).
\end{align*}
If $R>0$ is sufficiently large, so that 
$\tilde\la +C_n^2R^{-2}< \la(\Om\setminus K)$,
we obtain
\begin{equation*}
h_{a,V}(J_0u,J_0u)\le (\tilde\la+C_n^2R^{-2})(J_0u,J_0u)\le 
(\la(\Om\setminus K))(J_0u,J_0u), 
\quad u\in L.
\end{equation*}
(The inequality is strict if $J_0u\ne 0$.)
It follows from Corollary \ref{C:N-la-2}, that for the subspace 
$J_0(L)\subset C_0^\infty(\Om\cap B_{2R}(0))$ 
\begin{equation*}
\dim J_0(L)\le N\left(\la(\Om\setminus K); H_{a,V}|_{\Om\cap B_{2R}(0)}\right),
\end{equation*}
where $H_{a,V}|_{\Om\cap B_{2R}(0)}$ is the operator $H_{a,V}$ in $L^2(\Om\cap B_{2R}(0))$ 
with the Dirichlet boundary conditions, i.e. the operator defined by the closure 
of the quadratic form $h_{a,V}$ in $L^2(\Om\cap B_{2R}(0))$ from  $C_0^\infty(\Om\cap B_{2R}(0))$.  
Since $\Om\cap B_{2R}(0)$ is bounded, $H_{a,V}|_{\Om\cap B_{2R}(0)}$ has a discrete spectrum, 
so the number $N\left(\la(\Om\setminus K); H_{a,V}|_{\Om\cap B_{2R}(0)}\right)$
is finite. 

Now let us consider the multiplication-by-$J_0$ operator, restricted to $L$:
\begin{align*}\label{E:mult-by-J0}
M_{J_0}:L&\longrightarrow L^2(\Om\cap B_{2R}(0))\\
u&\longmapsto J_0u
\notag
\end{align*}
Then $\Ker M_{J_0}$ consists of functions from $L$ which vanish on $\supp J_0$
(in particular, in a neighborhood of $K$), so 
$\Ker M_{J_0}\subset C_0^\infty(\Om\setminus K)$.  Therefore, for any function
$u\in \Ker M_{J_0}$ both estimates \eqref{E:haV-on-L} and \eqref{E:haV-out} should be
satisfied. Since $\tilde\la<\la(\Om\setminus K)$, this is only possible if $u=0$.

So we conclude that $\Ker M_{J_0}=\{0\}$, so the map $M_{J_0}$ is injective. 
It follows that 
\begin{equation*}
N=\dim L=\dim J_0(L)\le N\left(\la(\Om\setminus K); H_{a,V}|_{\Om\cap B_{2R}(0)}\right).
\end{equation*}
This contradicts to the  assumption \eqref{E:la-3} which implies that $N$
can be arbitrarily large. Hence we proved \eqref{E:la-infty-below}. $\square$ 

\section{Two-sided estimates for the bottom of the essential spectrum}\label{S:2-sided-ess}

Now we turn to two-sided estimates of the bottom of the essential spectrum
for $H_{a,V}$ in $L^2(\Om)$. 
Denote $\Om_R=\Om\setminus \bB_R(0)$ and 
\begin{equation*} 
D_R=D_R(\Om,\ga,a,V)=D(\Om_R,\ga,a,V).
\end{equation*}
(cf. \eqref{E:D}). Since $D_R$ decreases as $R$ increases, we can define
\begin{equation}\label{E:D-infty}
D_\infty=D_\infty(\Om,\ga,a,V)=\lim_{R\to\infty} D_R(\Om,\ga,a,V).
\end{equation}
So $D_{\infty}\in [0,+\infty]$, and both $0$ and $+\infty$ can occur.
Clearly, $D_R$, $D_{\infty}$  increase with $\ga$ provided $\Om, a, V$
are fixed. Similarly, they increase with $\Om$, i.e. if $\Om\subset\Om'$,
then $D_R(\Om,\ga,a,V)\le  D_R(\Om',\ga,a,V)$,
and the same is true for $D_{\infty}$.

\begin{theorem}\label{T:2-sided-ess} 
There exists $C=C(\ga,n)>0$ such that
\begin{equation}\label{E:2-sided-ess}
C^{-1}D_{\infty}^{-2}\le \la_{\infty}\le C D_{\infty}^{-2},
\end{equation}
where $\la_{\infty}$ is the bottom of the essential spectrum of $H_{a,V}$
in $L^2(\Om)$.
\end{theorem}

{\bf Proof.} The result immediately folows from \eqref{E:D-infty} 
and Theorem \ref{T:2-sided}. $\square$

\begin{remark}
{\rm
Theorem \ref{T:discr2} with $\ga=const$ follows from Theorem \ref{T:2-sided-ess} because
the discreteness of spectrum is equivalent to the equality $\la_{\infty}=+\infty$, and 
the equality $D_{\infty}=0$ is equivalent to the corresponding conditions in 
Theorem \ref{T:discr2}.
}
\end{remark}

\section{A special class of operators}\label{S:example}

In this Section we will consider special magnetic Schr\"odinger operators
 $H_{a,V}$ in $L^2(\R^n)$, with the potentials of the form
 \begin{equation*}\label{E:special}
 a=a(x')=(0,\dots, 0, a_n(x')), \quad V=V(x'),
 \end{equation*}
 where $x'=(x^1,\dots, x^{n-1})$. In particular, $a$ and $V$ do not depend on the
 last coordinate $x^n$.  So the operator has the form
 \begin{equation}\label{E:HaV-special}
 H_{a,V} = -\De_{x'}  + \left(\frac{1}{i}\frac{\pa}{\pa x^n} + a_n(x')\right)^2 + V(x'), 
 \end{equation}
 We assume that  the  local regularity conditions 
 $V\in L^1_{loc}(\R^{n-1})$, $a_n\in L^2_{loc}(\R^{n-1})$
are satisfied, and, as above, $V\ge 0$, hence 
the self-adjoint operator $H_{a,V}$
is well defined through the quadratic form,
 
 Note that the magnetic field $B=da$ does not generally vanish 
 for such a potential $a$, but the corresponding skew-symmetric 
 matrix $(B_{jk})_{j,k=1}^n$ has a special form, with $B_{jk}=0$
 if $1\le j,k\le n-1$.
 
 Since the operator $H_{a,V}$ is invariant with respect to translations along the $x^n$ axis,
 making Fourier transform from $x^n$ to $\mu\in\R$, we obtain that $H_{a,V}$ is unitary
 equivalent to the following direct integral of self-adjoint operators
 \begin{equation}\label{E:tilde-H}
 \tilde{H}_{a,V}=\int^\oplus_{\R} H_{a,V}(\mu) \frac{d\mu}{2\pi}\;,
 \end{equation}
 where
 \begin{equation*}\label{E:HaV-mu}
 H_{a,V}(\mu) =  -\De_{x'}  + \left(\mu + a_n(x')\right)^2 + V(x'),
 \end{equation*}
 which is a Schr\"odinger operator without magnetic field in $L^2(\R^{n-1})$, 
 with a positive scalar
 potential
 \begin{equation*}
 V_\mu(x') =  \left(\mu + a_n(x')\right)^2 + V(x')
 =\mu^2 + 2\mu a_n(x') + a_n(x')^2 + V(x') .
 \end{equation*} 
 depending quadratically upon a parameter $\mu\in \R$. 
 Since the term $2\mu a_n(x')$ is dominated by the sum of the other terms in
 the right hand side, the regular perturbation theory applies, so it follows, in particular,
 that the bottom of the spectrum of $H_{a,V}(\mu)$ is a continuous function of $\mu$.
 Moreover, the bounded parts of whole spectrum of $H_{a,V}(\mu)$ are continuous
 with respect to $\mu\in\R$ in a natural sense. Due to the direct integral decomposition,
we see that all the spectrum of  $H_{a,V}$  is essential (i.e., it has no isolated points 
of finite multiplicity). 

Now let $\la, \la_\mu$ denote the bottoms of the spectra of the operators
$H_{a,V}$ and $H_{a,V}(\mu)$ respectively. Then it follows from 
the direct integral decomposition \eqref{E:tilde-H} 
and the arguments above, that
\begin{equation}\label{E:la-lamu}
\la=\inf \{\la_\mu|\;\mu\in\R\}.
\end{equation} 

Finally, similar to \eqref{E:D}, let us define
\begin{equation*}\label{E:tilde-D}
\tilde D=\tilde D(\ga, a,V)=\inf_{\mu\in\R} \sup_{Q_d}\left\{d:\; d^{n-2}\ge \inf_{F}\int_{Q_d\setminus F}
V_\mu(x')dx'\right\},
\end{equation*}
where $0<\ga< 1$ and the second infimum is taken over $F\subset Q_d$, satisfying
the negligibility condition $\capa(F)\le \ga\capa(Q_d)$. 
 
Now using the simplest version of Theorem \ref{T:2-sided} 
(without magnetic field, i.e. with $a\equiv 0$),
we immediately obtain
\begin{proposition}\label{P:2-sided} 
For every $\ga\in(0,1)$ there exists $C=C(\ga,n)>0$, such that
\begin{equation*}\label{E:2-sided-special}
C^{-1}\tilde D^{-2} \le \la \le C\tilde D^{-2}.
\end{equation*}
\end{proposition}

\end{document}